\documentclass[11pt]{article}
\usepackage{graphicx}

\def\be{\begin{equation}}
\def\ee{\end{equation}}
\def\ff#1{\mbox{\boldmath $#1$} }

\def\e{\epsilon}
\def\x{\ff{x}}

\begin{document}

\title{Mathematical Modelling and Parameter Optimization of Pulsating Heat Pipes}
\author{Xin-She Yang, Mehmet Karamanoglu, \\ 
School of Science and Technology, Middlesex University, London NW4 4BT, UK.
\and Tao Luan \\
School of Energy and Power Engineering, Shandong University, Jinan, China. \\
\and Slawomir Koziel \\
School of Science and Engineering, Reykjavik University, IS-103 Reykjavik, Iceland.
}

\date{}

\maketitle

\begin{abstract}
Proper heat transfer management is important to key electronic components in microelectronic
applications. Pulsating heat pipes (PHP) can be an efficient solution to such heat transfer problems.
However, mathematical modelling of a PHP system is still very challenging, due to the complexity and multiphysics nature of the system.
In this work, we present a simplified, two-phase heat transfer model, and our analysis shows that it can make good predictions
about startup characteristics. Furthermore, by considering parameter estimation as a nonlinear constrained optimization
problem, we have used the firefly algorithm to find parameter estimates efficiently.
We have also demonstrated that it is possible to obtain good estimates of key parameters using very limited
experimental data.
\end{abstract}

\noindent {\bf Citation detail:}
X. S. Yang, M. Karamanoglu, T. Luan, S. Koziel,
Mathematical Modelling and Parameter Optimization of Pulsating Heat Pipes,
{\it J. of Computational Science}, vol. 5, no. 2, pp. 119-125 (2014).

\newpage

\section{Introduction}

A pulsating heat pipe is essentially a small pipe filled with both liquid and vapour
and the internal diameter of the heat pipe is at the capillary scale \cite{Akac,Akac2}.
The liquids in the pipe can form segments or plugs, between vapour segments. When encountering heat,
part of the liquid many evaporate and absorb some heat, thus causing a differential pressure and driving
the movement of the plugs. When vapour bubbles meet a cold region, some of the vapour
may condensate, and thus releasing some heat. The loop can be open or closed, depending on the
type of applications and design. This continuous loop and process will
form an efficient cooling system if designed and managed properly for a given task.
Therefore, such systems have been applied to many applications in heat exchanger,
space applications and electronics, and they can potentially have even wider applications \cite{Zuo}.
On the other hand, the emergence of nanotechnology and the steady increase of the density of the
large-scale integrated circuits have attracted strong interests
in modelling heat transfer at very small scales, and the heat management of
microdevices has become increasingly important for next generation electronics
and miniaturization.

Both loop heat pipes (LHP) and pulsating heat pipes (PHP) may provide a promising
solution to such challenging problems, and thus have attracted renewed attention
in recent years \cite{Akac,Cheng,Fagh,Khan,Luan,Swan,Swan2,Zhir}. In many microelectronic applications,
conventional solutions to heat management problems often use fans, heat exchangers and
even water cooling. For examples, the fan-driven air circulation system in desktop computers
and many laptops have many drawbacks such as bulky sizes and potential failure of mechanical, moving parts.
In contrast, loop heat pipes have no mechanical driving system, and heat circulation
is carried out through the pipes, and thus such LHP systems can be very robust and long-lasting.
In addition, miniaturization and high-performance heat pipes systems are being developed \cite{Groll,McGlen,Reay,Yeom,Vasi,Zuo}.
Simulation tools and multiphase models have been investigated \cite{Rah,Naphon}.
All these studies suggested that LHP systems can have some advantages over traditional cooling systems.

A PHP system may often look seemingly simple; however, its working mechanisms are relatively
complex, as such systems involve multiphysics processes such as thermo-hydrodynamics, two-phase flow,
capillary actions, phase changes and others. Therefore, many challenging
issues still remain unsatisfactorily modelled. There are quite a few attempts in the literature to model
a PHP system with various degrees of approximation and success.

In this paper, we will use a mathematical
model based on one of the best models \cite{Swan,Swan2}, and will carry out some mathematical
analysis and highlight the key issues in the state of the art models.
In this paper, we intend to achieve two goals: to present a simplified mathematical model
which can reproduce most characteristics of known physics, and to provide a framework for
estimating key parameters from a limited number of measurements. The rest of
the paper is organized as follows: we first discuss briefly design optimization
and metaheuristic algorithms such as firefly algorithm, we then outline the main multiphysics processes in the
mathematical models. We then solve the simplified model numerically and compared with experimental
data drawn from the literature. After such theoretical analysis, we then form the parameter estimation as
an optimization problem and solve it using the efficient firefly algorithm for inversely estimating key parameters in a PHP system.
Finally, we highlight the key issues and discuss possible directions for further research.

\section{Design Optimization of Heat Pipes}
The proper design of pulsating heat pipes is important, so that the heat
in the system of concern can be transferred most efficiently. This also helps
to produce designs that use the least amount of materials and thus cost much less,
while lasting longer without the deterioration in performance.
Such design tasks are very challenging,
practical designs tend to be empirical and improvements tend to be incremental.
In order to produce better design options, we have to use efficient
design tools for solving such complex design optimization problems.
Therefore, metaheuristic algorithms are often needed to deal with such problems.

\subsection{Metaheuristics}

Metaheuristic algorithms such as the firefly algorithm and bat algorithm
are often nature-inspired \cite{Yang,YangBA},
and they are now among the most widely used algorithms for optimization. They have many advantages
over conventional algorithms, such as simplicity, flexibility, quick convergence
and capability of dealing with a diverse range of optimization problems.
There are a few recent reviews which are solely
dedicated to metaheuristic algorithms \cite{Yang,Yang2,YangDeb}.
Metaheuristic algorithms are very diverse, including
genetic algorithms, simulated annealing, differential evolution,
ant and bee algorithms, particle swarm optimization, harmony search, firefly algorithm,
cuckoo search, flower algorithm and others \cite{Yang2,Yang3,YangES,YangFPA}.

In the context of heat pipe designs, we have seen a lot of interests in the
literature \cite{Dob,Fagh}. However, for a given response, to identify the right
parameters can be considered as an inverse problem as well as an optimization problem.
Only when we understand the right working ranges of key parameters, we can start to
design better heat-transfer systems. To our knowledge,
this is the first attempt to use metaheuristic algorithms
to identify key parameters for given responses. We will use the firefly algorithm
to achieve this goal.

\subsection{Firefly Algorithm}

Firefly algorithm (FA) was first developed by Xin-She Yang in 2008 \cite{Yang,Yang3}
and is based on the flashing patterns and behavior of fireflies.
In essence, FA uses the following three idealized rules:

\begin{enumerate}
\item Fireflies are unisex so that one firefly will be attracted to other fireflies regardless of their sex.

\item The attractiveness is proportional to the brightness and they both decrease as their distance increases. Thus for any two flashing fireflies, the less brighter one will move towards the brighter one. If there is no brighter one than a particular firefly, it will move randomly in the form
    of a random walk.

\item The brightness of a firefly is determined by the landscape of the objective function.

\end{enumerate}

As a firefly's attractiveness is proportional to the light intensity seen by adjacent fireflies, we can now define the variation of attractiveness $\beta$ with the distance $r$ by
    \be \beta = \beta_0 e^{-\gamma r^2}, \ee
where $\beta_0$ is the attractiveness at $r=0$.

The movement of a firefly $i$ attracted to another more attractive (brighter) firefly $j$ is determined by
\be    \x_i^{t+1} =\x_i^t + \beta_0 e^{-\gamma r^2_{ij} } (\x_j^t-\x_i^t) + \alpha \; \ff{\e}_i^t, \ee
where the second term is due to the attraction. The third term is randomization with $\alpha$ being the randomization parameter, and $\ff{\e}_i^t$ is a vector of random numbers drawn from a Gaussian distribution or uniform distribution at time t. If $\beta_0=0$, it becomes a simple random walk. Furthermore, the randomization $\ff{\e}_i^t$ can easily be extended to other distributions such as L\'evy flights.

The L\'evy flight  essentially provides a random walk whose random step size $s$
is drawn from a L\'evy distribution
\be \textrm{L\'evy} \sim  s^{-\lambda}, \quad (1 < \lambda \le 3), \ee
which has an infinite variance with an infinite mean.
Here the steps essentially form
a random walk process with a power-law step-length distribution with a heavy tail.
Some of the new solutions should be generated by L\'evy walk around the
best solution obtained so far; this will speed up the local search. L\'evy flights
are more efficient than standard random walks \cite{Yang2}.

Firefly algorithm has attracted much attention \cite{Apo,Sayadi,Yang11}.
A discrete version of FA can efficiently solve non-deterministic polynomial-time hard, or NP-hard,
scheduling problems \cite{Sayadi}, while a detailed analysis has demonstrated the efficiency of FA for
a wide range of test problems, including multobjective load dispatch problems \cite{Apo,YangMOFA}.
Highly nonlinear and non-convex global optimization problems can be solved using firefly algorithm efficiently \cite{Gandomi,Yang11}.
The literature of firefly algorithms have expanded significantly, and Fister et al. provided
a comprehensive literature review \cite{Fister}.

\section{Mathematical Model for a PHP System}

\subsection{Governing Equations}

Mathematical modelling of two-phase pulsating flow inside a pulsating heat pipe involves many processes, including
interfacial mass transfer, capillary force, wall shear stress due to viscous action, contact angles,
phase changes such as evaporation and condensation, surface tension, gravity and adiabatic process.
All these will involve some constitutive laws and they will be coupled with fundamental laws
of the conservation of mass, momentum and energy, and thus resulting in a nonlinear system of
highly coupled partial differential equations. Consequently, such a complex model can
lead to complex behaviour, including nonlinear oscillations and even chaotic characteristics \cite{Maez,Qu,Saku,Swan,Wong,Zhang}.

A mathematical model can have different levels of complexity, and often a simple model
can provide significant insight into the working mechanism of the system and its
behaviour if the model is constructed correctly with realistic conditions. In most cases,
full mathematical analysis is not possible, we can only focus on some aspects of the
model and gain some insight into the system.

A key relationship concerning the inner diameter of a typical pulsating heat pipe is the range of capillary length,
and many design options suggest to use a diameter near the critical diameter \cite{Dob}
\be d_c = 2 \sqrt{\frac{\sigma}{g (\rho_l-\rho_v)}}, \ee
where $\sigma$ is the surface tension (N/m), while $g$ is the acceleration due to gravity,
which can be taken as $9.8$ m/$s^2$. $\rho_l$ and $\rho_v$ are the densities of liquid and vapour, respectively.
In the rest of the paper, we will focus on one model which may be claimed as the state-of-the-art, as it is based on
the latest models \cite{Swan,Swan2}, with a simplified model for a system of
liquid plugs and vapour bubbles as described in \cite{YangICCS2012}.

\subsubsection{Temperature Evolution}

The temperature $T_v$ in a vapour bubble is governed by the energy balance equation
\be m_v c_{vv} \frac{d T_v}{dt}=-h_{lfv} (T_v -\tau) L \pi (d_i -2 \delta)-r_m h_v L \pi (d_i - 2 \delta) -p_v \frac{d V}{dt}, \ee
where
\be r_m = \frac{2 \sigma_0}{(2- \sigma_0)} \frac{1}{\sqrt{2 \pi R}} (\frac{p_v}{\sqrt{T_v}}-\frac{p_l}{\sqrt{\tau}}). \ee
Here $L$ is the mean plug length, $d_i$ is the inner diameter of the pipe, and $\delta$ is the thickness of the thin liquid film.
$h$ corresponds to the heat transfer coefficient and/or enthalpy in different terms. Here $\sigma_0$ is a coefficient and $R$ is the gas constant.

In addition, the volume in the above equation is given by
\be V=\frac{\pi d_i^2}{4} (L+L_v), \ee
where the length of the vapour bubble is typically $L_v=0.02$ m, and $\sigma \approx 1$.

The temperature $\tau$ in the liquid film is governed by
\be m_{f} c_{vl} \frac{d \tau}{dt} =- h_{lfw} (\tau-T_w) L \pi d_i + h_{lfv} (T_v-\tau) L \pi (d_i - 2 \delta) + r_m h_v L \pi (d_i - 2 \delta), \ee
where $T_w$ is the initial wall temperature.

\subsubsection{Mass Transfer}
For a vapour bubble, rate of change in mass is governed by the conservation of mass
\be \frac{d m_v}{dt} = - \pi (d_i - 2 \delta) r_m. \ee
As $m_v=\rho_v V$, we have
\be \rho_v \frac{dV}{dt} = \pi (d_i-2 \delta) r_m - r_v \pi d_i L_v, \ee
where $\rho_v \approx 1 $ Kg/m$^3$ and the transfer rate $r_v=0$ if $T_w>T_v$.

\subsubsection{Motion of a Plug}
The position $x_p$ of the liquid plug is governed by the momentum equation
\be m_p \frac{d^2 x_p}{d t^2} = \frac{\pi}{4} (d_i- 2 \delta)^2 (p_{v1} -p_{v2}) - \pi d_i L_p s_w + m_p g, \ee
where
\be m_p =\rho_l \frac{\pi d_i^2}{4} (L_0-x_p), \ee
with $L_0 \approx 25 d_i$ and $\rho_l=1000$ Kg/m$^3$. Here the term $p_{v1}-p_{v2}$ is the differential
pressure between the two sides of the plug.

In the above equations, we also assume the vapour acts as an ideal gas, and we have
\be p_v=\frac{m_v R T_v}{V}. \ee
The shear stress $s_w$ between the wall and the liquid plug is given by
an empirical relationship
\be s_w= \frac{1}{2} C_f \rho_l v_p^2, \ee
where $C_f\approx 16/Re$ when $Re=\rho_l v_p d_i/\mu_l \le 1180$; otherwise,
\be C_f \approx 0.078 Re^{-1/4}. \ee
Here, we also used that $v_p=d x_p/dt$.

The mathematical formulation can include many assumptions and simplifications.
When writing down the above equations, we have tried to incorporate most the
constitutive relationships such contact angles, capillary pressures and viscous force into the equations directly so that we can
have as fewer equations as possible. For details, readers can refer to \cite{Swan,Swan2,Dob}. This way, we can focus on
a few key equations such as the motion of a plug and temperature
variations, which makes it more convenient for later mathematical analysis.

\subsection{Typical Parameters} \label{section-100}

The properties of the fluids and gas can be measured directly, and typical values can be
summarized here, which can be relevant to this study. Most of these values have been drawn from earlier
studies \cite{Dob2,Swan,Yaws}.

Typically, we have $L=0.18$ m, $d_i=3.3 \times 10^{-3}$ m,
$\delta=2.5 \times 10^{-5}$ m. The initial temperature is $T_v \approx \tau=20^{\circ}$ C,
while the initial wall temperature is $T_w =40^{\circ}$ C. The initial pressure
can be taken as $p_v \approx 5.5816 \times 10^4$ Pa. $h_{lfw}=1000 $W/m$^2$K,
$h_{v}=10$ W/m$^2$K. $c_{vl}=1900 $J/Kg ${}^{\circ}$C, $R=8.31$
and $c_{vv}=1800 $J/Kg ${}^{\circ}$C.

In addition, we have the initial values: $m_{v0}= p_{v0} V_0/(R_v (T_{v0}+273.15)$ where
$T_{v0}=20$ ${}^{\circ}$C, $p_{v0}=10^{5}$ Pa, $R_v=461$ J/Kg K,
and $m_{f0} \approx m_{v0}/10$.

\section{Nondimensionalization, Analysis and Simulation}

After some straightforward mathematical simplifications,
the full mathematical model can be written as the following equations:
\be
\begin{array}{rcl}
 m_v c_{vv} \frac{d T_v}{dt}&=&-h_{lfv} (T_v -\tau) L \pi (d_i -2 \delta)-r_m h_v L \pi (d_i - 2 \delta) \\
 & &  -p_v \frac{\pi}{\rho_v}[(d_i-2 \delta) r_m - d_i L_v r_v], \\ \\
 m_{f} c_{vl} \frac{d \tau}{dt} &=&- h_{lfw} (\tau-T_w) L \pi d_i + h_{lfv} (T_v-\tau) L \pi (d_i - 2 \delta) \\
 & & + r_m h_v L \pi (d_i - 2 \delta), \\ \\
\frac{d m_v}{dt} &=& - \pi (d_i - 2 \delta) r_m, \\ \\
m_p \frac{d^2 x_p}{d t^2} &=& \frac{\pi}{4} (d_i- 2 \delta)^2 (p_{vk} -p_{vk}) - \pi d_i L_p s_w + m_p g. \\
\end{array}
\label{math-100}
\ee
Though analytically intractable, we can still solve this full model numerically using
any suitable numerical methods such as finite difference methods, and see how the system behaves under various conditions with various values of parameters.
Preliminary studies exist in the literature \cite{Swan,Dob2}, and we have tested our system that it can indeed reproduce these results using the typical
parameters given in Section \ref{section-100} \cite{YangICCS2012}. For example, the mean temperature in a plug is consistent with the results in \cite{Zhang}, while the oscillatory behaviour is similar to the results obtained by others \cite{Khan,Dob,Zhang,Swan2,Yuan}.

\subsection{Non-dimensionalized Model}

For the current purpose of parameter estimation,  the system can be considered stationary under appropriate conditions,
and thus many terms can be
taken as constants. Then, the full mathematical model (\ref{math-100}) can be non-dimensionalized and written as
\be
\begin{array}{rcl}
u \frac{dT_v}{dt} &=& a (T_v-\tau) -\alpha_1 -\alpha_2, \\ \\
\frac{d\tau}{dt} &=& b (T_v-\tau) - \epsilon \tau + \alpha_3, \\  \\
\frac{d u}{dt} &=& - \Delta, \\ \\
\frac{d^2 x_p}{dt^2} &=& \frac{1}{(\beta_1-\beta_2 x_p)}\Big[\beta -\gamma (\frac{d x_p}{dt})^2\Big],
\end{array}
\ee
where
\be a=-\frac{h_{lfv} L \pi (d_i - 2 \delta)}{ c_{vv}}, \quad \alpha_1=\frac{r_m h_v L \pi (d_i -2 \delta)}{ c_{vv}}, \ee
\be \alpha_2=\frac{p_v \pi}{\rho_v  c_{vv}} [(d-2 \delta) r_m - d_i L_v r_v], \quad
 b=-\frac{h_{lfv} L \pi (d_i -2\delta)}{m_f c_{vl}}, \ee
 \be  \epsilon=\frac{h_{lfw} L \pi d_i}{m_f c_{vl}}, \quad \alpha_3=\frac{r_m h_v L \pi (d_i-2 \delta)}{m_f c_{vl}},  \ee
\be \Delta=\frac{\pi (d_i-2 \delta ) r_m}{m_{v0}}, \quad \beta=\frac{\pi (d_i- 2\delta)^2 (p_{v1}-p_{v2})}{4} +g, \ee
\be \gamma=\frac{\pi d_i L_p C_v \rho_l}{2}, \quad \beta_1=\frac{\rho_l L_0 \pi d_i^2}{4}, \quad \beta_2= \frac{\rho_l \pi d_i^2}{4}. \ee
In the above formulation, $u$ is the dimensionless form of $m_v$, that is $u=m_v/m_{v0}$.

\subsection{Simplified Model}
\begin{figure}[h]
\centerline{\includegraphics[width=4in,height=3.15in]{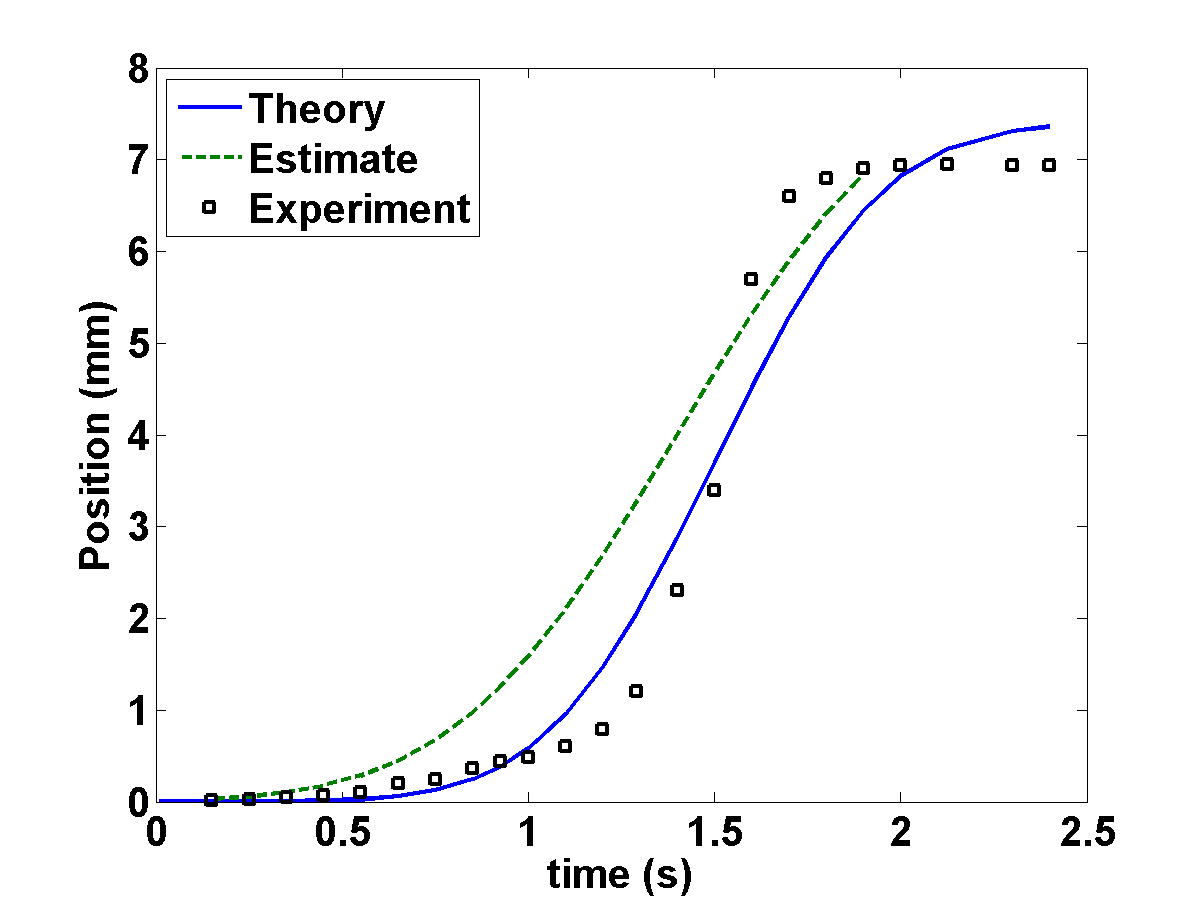}}
\caption{Comparison of the full numerical results with approximations and experimental data \cite{Swan2}. \label{ykfig1} }
\end{figure}

Under the assumptions of stationary conditions and $x_p/L_0 \ll 1$ (for example at the startup stage),
the last equation is essentially decoupled from the other equations.  In this case, we have an approximate
but simplified model
\be
\begin{array}{ll}
\frac{d\tau}{dt} = b (T_v-\tau) - \epsilon \tau + \alpha_3, \\  \\
\frac{d^2 x_p}{dt^2} = \frac{1}{\beta_1}\Big[\beta -\gamma (\frac{d x_p}{dt})^2\Big],
\end{array}
\ee
Now the second equation can be rewritten as
\be \frac{d^2 x_p}{d t^2}=A - B (\frac{d x_p}{dt})^2, \quad  A=\frac{\beta}{\beta_1}, \quad B=\frac{\gamma}{\beta_1}. \ee
For the natural condition $x_p=0$ at $t=0$, the above equation has a short-time solution
\be x_p =\frac{1}{B} \ln \Big[\cosh(\sqrt{A B} \; t) \Big], \ee which
is only valid when $x_p/L_0$ is small and/or $t$ is small. Indeed, it gives some main characteristics of the
startup behaviour as shown as the dashed curve in Fig. \ref{ykfig1}.

Furthermore, the first equation for the liquid plug temperature is also decoupled from the PDE system. We have
\be \frac{d \tau}{dt} =Q_1 - Q_2 \tau, \quad Q_1= b T_v + \alpha_3, \quad Q_2 =b+ \epsilon. \ee
For an initial condition $\tau=0$ at $t=0$, we have the following solution
\be \tau=Q_1 (1-e^{-Q_2 t}), \ee
which provides a similar startup behaviour for small times as that in \cite{Zhang}.

The approximation solutions are compared with the full numerical solutions and experimental data points, and
they are all shown in Fig.~\ref{ykfig1}.
Here the experimental data points were based on \cite{Swan,Swan2,Dob}. It is worth pointing that the approximation
can indeed provide some basic characteristics of the fundamental behaviour of the startup, and the full numerical results can
also give a good indication of the final steady-state value.

This simplified model also confirms that some earlier studies using mainly $x_p$ as the dependent variable
can indeed give good insight into the basic characteristics of the complex system. For example, Wong {\it et al.}
used an equivalent viscous damping system without considering the actual heat transfer system \cite{Wong}, and their
main governing equation can be written as
\be \frac{d^2 x}{d t^2} + a \frac{d x}{dt} + b (k + t) x =0, \ee
where $a, b, k$ are constants. This system has typical features of a viscous damping system. On the other hand,
Yuan {\it et al.} used a system with primarily a second-order nonlinear ODE for $x_p$
\be \frac{d^2 x_p}{d t^2} + \frac{2 C}{d_i} \Big(\frac{d x_p}{dt}\Big)^2 + \frac{2 g}{L} x_p=\frac{(p_{v1}-p_{v2})}{L \rho_l}, \ee
where $C$ is the friction coefficient. This system can have even richer characteristics of oscillatory behaviours \cite{Yuan}.

In comparison with earlier simpler models, our simplified model can produce even richer dynamic features of the
system without using much complex mathematical analysis, thus provides a basis for further realistic analysis and simulations.

\section{Inverse Parameter Identification and Optimization}

The aim of an inverse problem is typically  to estimate important parameters of structures and materials, given
 observed data which are often incomplete. The target is to minimize the differences between
  observations and predictions, which is in fact an optimization problem.
  To improve the quality of the estimates, we have to combine a wide range of known information,
  including any prior knowledge of the structures, available data. To incorporate all useful
  information and carry out the minimization, we have to deal with a multi-objective optimization
  problem. In the simplest case, we have to deal with a nonlinear least-squares problem with
  complex constraints \cite{Samb,Samb2}.

The constraints for inversion can include the realistic ranges of physical parameters,
 geometry of the structures, and others. The constraints are typically nonlinear,
 and can be implicit or even black-box functions.
 Under various complex constraints, we have to deal with a nonlinear, constrained, global optimization problem.
 In principle, we can then solve the formulated constrained problem by any efficient optimization
 techniques \cite{YangICCS2012,Yang2}. However, as the number of the degrees of freedom in inversion is typically large,
 data are incomplete, and non-unique solutions or multiple solutions may exist;
 therefore, meteheuristic algorithms are particularly suitable.

\subsection{Parameter Estimation as an Inverse Problem}

Inverse problems in heat transfer and other disciplines tend to find the best parameters of
 interest so as to minimize the differences between the predicted results and observations.
 In the simplest case, a linear inverse problem can be written as

\be \textbf{u}=\textbf{Kq} + \textbf{w}, \ee
 where \textbf{u} is the observed data with some noise \textbf{w}, \textbf{q} is the
 parameter to be estimated, and \textbf{K} is a linear operator or mapping, often
 written as a matrix. An optimal solution \textbf{q*} should minimize the
 norm $||\textbf{u}-\textbf{Kq}||$. In case when \textbf{K} is rectangular, we have the best estimate

 \be \textbf{q}^*=(\textbf{K}^{T} \textbf{K})^{-1} \textbf{K}^{T} \textbf{u}. \ee

 But in most cases, \textbf{K} is not invertible, thus we have to use other techniques
 such as regularization.  In reality, many inverse problems are nonlinear and data are
 incomplete. Mathematically speaking,  for a domain $\Omega$ such as a structure with
 unknown but true parameters \textbf{q}* (a vector of multiple parameters),  the
 aim is find the solution vector \textbf{q} so that the predicted values y${}_{i}$ (i=1,2, \dots , n),
 based on a mathematical model $y=\phi(\x,p)$ for all $\x \in \Omega$,  are close to the observed values ${\bf d}=(d_1, d_2, ..., d_n)$ as possible.

 The above inverse problem can be equivalently written as
 a generalized least-squares problem \cite{Samb,YangICCS2012}

\be   \min f=|| d-\phi(\x,\textbf{q})||{}^{2}, \ee
or
\be     \min \quad \sum _{i} [d_{i}^{} -\phi (x_{i} ,q)]^{2}, \label{yk-leastsq}\ee

 Obviously, this minimization problem is often subjected to a set of complex constraints.
 For example, physical parameters must be within certain limits. Other
 physical and geometrical limits can also be written as nonlinear constraints.
 In general, this is equivalent to the following general nonlinear constrained optimization problem
\be \min \quad f(\x,{\bf q},{\bf d}^{} )  \label{yk-kopt}\ee
subject to
\be h_{j} (\x,{\bf q})=0,\quad (j=1,...,J),\quad g_{k} (\x,{\bf q})\le 0,\quad (k=1,...,K). \ee
where J and K are the numbers of equality and inequality constraints, respectively.
Here all the functions $f$, $h_{j}$, and $g_{k}$ can be nonlinear functions. This is a nonlinear, global optimization for \textbf{q}. In the
case when the observations data are incomplete, thus the system is under-determined,
some regularization techniques such as Tikhonov regularization are needed.
Therefore, the main task now is to find an optimal solution to approximate the true parameter set \textbf{q}*.

From the solution point of view,  such optimization can be in principle solved using any efficient optimization
 algorithm. However, as the number of free parameters tends to be very large,
 and as the problem is nonlinear and possible multimodal, conventional algorithms
 such as hill-climbing usually do not work well. More sophisticated metaheuristic
 algorithms have the potential to provide better solution strategies.

Though above inverse problems may be solved if they are well-posed,
however, for most inverse problems, there are many challenging issues.
First, data are often incomplete, which  leads to non-unique solutions; consequently,
the solution techniques are often problem-specific, such as Tikhonov
regularization. Secondly, inverse problems are highly nonlinear
and multimodal, and thus very difficult to solve. In addition, problems are often
large-scale with millions of degrees of freedom, and thus requires efficient
algorithms and techniques. Finally, many problems are NP-hard, and there is no
efficient algorithms of polynomial time exists. In many of these cases, metaheuristic algorithms
such as genetic algorithms and firefly algorithm could be the only alternative. In fact, metaheuristic algorithms
are increasingly popular and powerful \cite{Yang,YangICCS2012}. In the rest of this paper,
we will use the firefly algorithm discussed earlier to carry out the parameter estimations
for pulsating heat pipes.

\subsection{Least-Squares Estimation}

Now we try to solve the following problem for parameter estimation: Suppose we measure the response (i.e., the location $x_p$ of a plug)
of a pulsating heat pipe, as the results presented in Fig.~\ref{ykfig1},
can we estimate some key parameters using these results through our simplified mathematical model?
One way to deal with this problem is to consider it as a nonlinear, least-squares best-fit problem as described earlier
in Eq. (\ref{yk-leastsq}).

Let us focus on the key parameters $L$, $T_v$, $T_w$, $d_i$, and $p_v$.
When we carry out the estimation using the least squares methods, we can only get
crude estimates as shown in Table 1.
\begin{table}
\begin{center}
\caption{Least-Squares Parameter Estimation.}
\begin{tabular}{c|c|c}
\hline
Parameters & Estimates & True values \\ \hline
L & 0.15 $\sim$ 0. 24 & 0.18  (m)\\
$d_i$ & 0.001 $\sim$ 0.005 & 0.0033  (m) \\
$T_v$ & 10 $\sim$ 27 & 20 (${}^{\circ}$C)\\
$T_w$ & 20 $\sim$ 49 & 40 (${}^{\circ}$C) \\
$p_v$ & 40 $\sim$ 129 & 100 (kPa) \\ \hline
\end{tabular}
\end{center}
\end{table}
Here we have used 25 data points to establish estimates for 5 key parameters.
The wide variations of these parameters, though near the true values, suggest that there are
insufficient conditions for the inverse problem to have unique solutions. This can be
attributed to the factors that the measured data points were sparse and incomplete, no
constraints were imposed on the ranges of the parameters, and additional constraints might be needed to
make the inverse problem well-posed.

\subsection{Constrained Optimization by the Firefly Algorithm}

In order to get better and unique solutions to the problem, to increase of the number
of measurements is not the best choice, as experiments can be expensive and time-consuming.
Even with the best and most accurate results, we may still be unable to get unique solutions.
We have to make this problem well-posed by imposing enough conditions.

First, we have to impose more stringent bounds/limits. Thus, we apply
$L \in [0.15, 0.22]$ m, $d_i \in [0.002, 0.004]$ m,
$T_v \in [15, 25]$ ${}^{\circ}$C, $T_w \in [35, 45]$
${}^{\circ}$C, and $p_v \in [80, 120]$ kPa.
Then, we have to minimize the parameter variations, in addition to the residual sum of squares. Thus, we have
\be     \textrm{Minimize } \quad \sum _{i} \Big[d_{i}^{} -\phi (x_{i} ,q)\Big]^{2}+ \sum_{k=1}^5 \sigma^2_k, \label{yk-leastsqnew}\ee
where $\sigma_k^2$ is the variance of parameter $k$. Here we have 5 parameters
to be estimated. Now we have a constrained optimization problem with imposed simple bounds.

By using the firefly algorithm with $n=20$, $\gamma=1$, $\beta=1$ and $5000$ iterations,
we can solve the above constrained optimization. We have run the results 40 times so as to
obtained meaningful statistics, using the same set of 25 observed data points, and the results are summarized in Table 2.

\begin{table}
\begin{center}
\caption{Parameter Estimation by Constrained Optimization.}
\begin{tabular}{c|c|c}
\hline
Parameters & Estimates &  True values \\
& (mean $\pm$ standard deviation) &  \\ \hline
L & 0.17 $\pm$ 0.1 & 0.18  (m)\\
$d_i$ & 0.0032 $\pm$ 0.0004 & 0.0033  (m) \\
$T_v$ & 20.4 $\pm$ 0.9 & 20 (${}^{\circ}$C)\\
$T_w$ & 39.2 $\pm$ 1.1 & 40 (${}^{\circ}$C) \\
$p_v$ & 98 $\pm$ 8.8 & 100 (kPa) \\ \hline
\end{tabular}
\end{center}
\end{table}
We can see clearly from the above results that unique parameter estimates can be
achieved if sufficient constraints can be imposed realistically,
and the good quality estimates can be obtained even for a number of sparse measurements.
This implies that it is feasible to estimate key parameters of a complex PHP system from experimental data
using the efficient optimization algorithm and the correct mathematical models.

\section{Discussions and Conclusions}

Though a pulsating heat pipe system can be hugely complex as it involves multiphysics processes,
we have shown that it is possible to formulate a simplified mathematical model under
one-dimensional configuration, and such a model can still be capable to reproduce the fundamental
characteristics of time-dependent startup and motion of liquid plugs in the heat pipe.
By using scaling variations and consequently non-dimensionalization, we have
analysed the key parameters and processes that control the main heat transfer process.
Asymptotic analysis has enabled us to simplify the model further concerning
small-time and long-time behavior, so that we can identify the factors that
affect the startup characteristics. Consequently, we can use the simplified model
to predict and then compare the locations of a plug with experimental data.
The good comparison suggests that a simplified model can work very well
for a complex PHP system.

Even though we have very good results, however, there are still some significant differences
between the full numerical model and experimental data, which can be attributed to unrealistic
parameter values, incomplete data, oversimplified approximations and unaccounted experimental settings.
Further work can focus on the extension of the current model to 2D or 3D configuration
with realistic geometry and boundary conditions.

By treating parameter estimation as an inverse problem and subsequently a nonlinear constrained
optimization problem, we have used the efficient metaheuristic algorithms such as firefly algorithm
to obtain estimates of some key parameters in a 1D PHP system. We have demonstrated that a least-squares approach is
not sufficient to obtain accurate results because the inverse problem was not well-posed, with insufficient constraints.
By imposing extra proper constraints on parameter bounds and also minimizing the
possible variations as well as minimizing the best-fit errors, we have obtained
far more accurate estimates for the same five key parameters, and these parameter estimates are comparable
with their true values.

There is no doubt that further improvement concerning the mathematical model and parameter estimation
will provide further insight into the actual behaviour of a PHP system, and subsequently
allows us to design better and more energy-efficient PHP systems.

\section*{References}
\bibliographystyle{elsarticle-num}

\end{document}